\documentstyle[a4,11pt,righttag,amssymb,verbatim]{amsart}

\newtheorem{Th}{Theorem}[section]

\newtheorem{Lem}[Th]{Lemma}
\newtheorem{Prop}[Th]{Proposition}
\newtheorem{Clm}[Th]{Claim}

\newtheorem{th}{Theorem} 
\newtheorem{prop}{Proposition}

\numberwithin{equation}{section}
\renewcommand{\theequation}{\thesection.\arabic{equation}}

\def\cB{{\cal B}} \def\cC{{\cal C}} 
  \def\Inn#1{\operatorname{Inn}(#1)}
\def\Fix{\operatorname{Fix}} \def\Cen{\operatorname{Cen}}

\def\f{\varphi} \def\s{\sigma} \def\eps{\varepsilon} 
\def\aut#1{\operatorname{Aut}(#1)}
\def\str#1{\langle #1\rangle} \def\rank#1{\operatorname{rank}#1}
\def\av#1{\overline{#1}}     \def\avst#1{\overline{\mathstrut #1}}
\def\inv{{}^{-1}} \def\id{\operatorname{id}}
\def\Z{\boldsymbol{\text Z}}

\begin{document}

\title[The automorphism tower of a free group]{The automorphism
tower of a free group} 
\thanks{Supported by Russian Foundation of Fundamental Research
Grant 96-01-00456}
\author{Vladimir Tolstykh}
\address{Department of mathematics \\ Kemerovo State University\\
Krasnaja, 6  \\ 650043 Kemerovo \\ Russia}
\date{September 19, 1997}
\maketitle

\begin{abstract}
We prove that the automorphism group of any non-abelian
free group $F$ is complete. The key technical step in the proof:
the set of all conjugations by powers of primitive elements
is first-order parameter-free definable in the group
$\aut F.$
\end{abstract}

\section*{\it Introduction}

In 1975 J.~Dyer and E.~Formanek \cite{DFo} had proved
that the automorphism group of a finitely generated
non-abelian free group $F$ is complete (that is, it is
centreless and all its automorphisms are inner) and
so $\aut{\aut F} \cong \aut F$. They noted that
their research was stimulated by G. Baumslag, who
conjectured that the automorphism tower of a finitely
generated free group is very short. New proofs for the
result of Dyer and Formanek were given in 1990 by
D.~G.~Khramtsov \cite{Khr} and E.~Formanek \cite{Fo}.

The objective of this paper is to generalize the result
of Dyer and Formanek from finitely generated non-abelian
free groups to arbitrary non-abelian free groups.

Let $F$ be a free non-abelian group.  We obtain a
group-theoretic characterization of conjugations by
powers of primitive elements in $\aut F.$  Our key
technical results can be summarized in model-theoretic
terms as follows: the set of all conjugations by
powers of primitive elements is first-order
parameter-free definable in the group $\aut F$
(Theorem \ref{DefOfConjsByPPE}). The latter means that
there is a first-order formula with one free variable
in the language of groups such that its realizations in
$\aut F$ are exactly conjugations just mentioned.
Therefore the subgroup of all conjugations (inner
automorphisms of $F$) is a characteristic subgroup of
$\aut F.$ This implies that the group $\aut F$ is
complete (Theorem \ref{MainRes}).

The main technical tool in the proof of Theorem
\ref{DefOfConjsByPPE} is the use of conjugacy classes
of involutions based on a characterization of
involutions in $\aut F$ given by  J.~Dyer and
G.~P.~Scott in \cite{DSc}.  An important role is
played by involutions of the following sort. Let $x$
be a primitive element of $F$ and $F=\str x * C$ a
free factorization of $F.$ Then an automorphism of $F$
which inverts $x$ and takes each element in $C$ to its
conjugate by $x,$ is an involution. We call any
involution obtained in such a way a {\it
quasi-conjugation}, since it acts as conjugation on  a
`large' subgroup of $F.$

Suppose $\rank F> 2.$ For given a quasi-conjugation $\varphi \in \aut F$
let $\Pi$ be the set of all automorphisms of $F$ of the form
$\s\s',$ where $\s$ and $\s'$ both commute with $\varphi$
and are conjugate. We prove (Proposition \ref{Conjs:rank=>3})
that if $\varphi$ is defined by $x$ and $C$ then conjugations
by powers of $x$ are exactly the members of the centralizer
of $\Pi$ in $\aut F$ which are not involutions. Similarly,
we characterize conjugations by powers of primitive
elements in the case when $\rank F=2$  (Proposition \ref{Conjs:rank=2}).
These results reduce the problem of first-order characterization
of conjugations by powers of primitive elements to a
characterization of quasi-conjugations. The latter
problem is solved in Section \ref{q-conjs}: we
characterize the class of all quasi-conjugations in
terms of products of conjugacy classes.

We say that a subset of a group is {\it
anti-commutative} if its elements are pairwise
non-commuting.  In the case when $\rank F > 2$ we
prove that the class of all quasi-conjugations is the
unique anti-commutative conjugacy class of
involutions such that, for every anti-commutative
conjugacy class $K'$ of involutions, all involutions in $KK'$ are
conjugate (Proposition \ref{DefOfQsConjs}). When
$\rank F=2$ the class of all quasi-conjugations is the
unique anti-commutative conjugacy class $K$ of
involutions such that elements in $K$ are not squares
(Proposition \ref{rank(F)=2}).  It enables us to
construct a first-order formula characterizing
quasi-conjugations in $\aut F,$ and hence to do the
same for conjugations by powers of primitive elements
in $F.$

The results on conjugacy classes of involutions in
$\aut F$ needed in the main body of the paper are
considered in Sections~\ref{invs} and~\ref{acc}.

The author is very grateful to his colleagues in Kemerovo
University Oleg Belegradek, Valery Mishkin, and Peter
Biryukov for reading of the first draft of this paper and
helpful comments.  The main result of the paper
(Theorem \ref{MainRes}) was announced in the
abstract~\cite{To}.

\section{\it Notation and preliminaries}

In what follows $F$ stands for a free
non-abelian group. The free abelian group $F/[F,F]$ of the
same rank is denoted by $A.$ The natural homomorphism $w
\mapsto \av w$ from the group $F$ to $A$ provides the
homomorphism $\aut F \to \aut A.$ To denote this
homomorphism we shall be using the same symbol
$\av{\phantom{a}}.$ We shall also say that an
automorphism $\varphi$ of $F$ {\it induces} the
automorphism $\av\varphi \in \aut A.$

We shall use the following fact.

\begin{Prop} \label{Specht?} \mbox{\rm (\cite{Spe})}
Let $G$ be a centreless group. Then the group $\aut G$ is complete
if and only if the subgroup $\Inn F$ of inner automorphisms
{\rm(}conjugations{\rm)} is characteristic in $\aut G.$
\end{Prop}

In this paper we prefer to call elements in the subgroup
$\Inn F$ `conjugations' rather than `inner automorphisms of $F$'.

It is convenient to formulate our main technical
results using a model-theoretic notion of definable
set (see \cite[ch. II]{Ho}). An $n$-ary relation $S$
on a group $G$ is said to be {\it first-order definable
without parameters} in $G$ (or, for short {\it parameter-free
definable}) if there is a first-order formula $\chi(v_1,\ldots,v_n)$
in the language of groups $\{\cdot,{}\inv,1\}$ such that
$S$ is the set of all $n$-tuples $(a_1,\ldots,a_n)$ in $G$
realizing $\chi(a_1,\ldots,a_n)$ in the group $G.$
For example, the centre of $G$ is parameter-free definable
by the formula $(\forall u)(vu=uv).$ Clearly, every parameter-free
definable relation on $G$ admits a description in terms
of the group operation.  Therefore

\begin{Prop} \label{Def=>Char}
Any parameter-free definable subset of a group $G$ is invariant
under all automorphisms of $G$ and hence generates
a characteristic subgroup of $G.$
\end{Prop}

We shall prove that for any non-abelian free group $F$
the set of all quasi-conjugations and the set of all conjugations
by powers of primitive elements both are parameter-free definable
in $\aut F$ and hence are invariant under automorphisms of
$\aut F.$

\section{\it Involutions} \label{invs}

In \cite{DSc} J.~Dyer and G.~P.~Scott obtained a description
of automorphisms of $F$ of prime order. For involutions
that description yields the following

\begin{Th} \mbox{\rm \cite[p. 199]{DSc} } \label{CanForm}
For every involution $\varphi$ in
the group $\aut F$ there is a basis $\cB$ of $F$ of the form
$$
\{u : u \in U \} \cup \{z, z' : z \in Z \} \cup
\{ x, y : x \in X, y \in Y_{x} \}
$$
on which $\varphi$ acts as follows
\begin{align} \label{eqCanForm}
&\quad \varphi u = u, \quad u \in U, \\
& \begin{cases}
\varphi z = z' & z \in Z, \\
\varphi z'= z, & \\
\end{cases}  \nonumber \\
&\begin{cases}
\varphi x =x \inv,   & x \in X,\\
\varphi y =x y x\inv, & y \in Y_x.
\end{cases} \nonumber
\end{align}
Specifically, the fixed point subgroup of $\varphi,$ $\Fix(\varphi),$ is the
subgroup $\str{u : u \in U},$ and hence is a
free factor of $F.$
\end{Th}

We shall call a basis of $F$ on which $\varphi$ acts
similar to (\ref{eqCanForm}) a {\it canonical} basis
for $\varphi.$ In view of (\ref{eqCanForm}) one can
partition every canonical basis $\cB$ for $\varphi$ as follows
\begin{equation} \label{eqBDivision}
\cB=U(\cB) \cup Z(\cB) \cup \{z' : z \in Z(\cB)\}\cup X(\cB) \cup
\bigcup_{x \in X(\cB)} Y_x(\cB).
\end{equation}
We shall also call any set of the form $\{x\} \cup Y_x,$
where $x \in X(\cB)$ a {\it block} of $\cB,$ and
the cardinal $|Y_x|+1$ the {\it size} of a block.
The subgroup generated by the set $Y_x$ will be
denoted by $C_x$ ($\varphi$ operates on this subgroup as
conjugation by $x$), and the subgroup generated by the block
$\{x\} \cup Y_x$ will be denoted by $H_x.$ Sometimes
we shall be using more `accurate' notation like
$C^{\varphi}_x$ or $H^{\varphi}_x.$ The set $U(\cB)$ will be
called the {\it fixed part} of $\cB.$

Clearly, if $\cB$ and $\cC$ are some canonical
bases for involutions $\varphi,\psi,$ respectively, and
the action of $\varphi$ on $\cB$ is isomorphic to the
action of $\psi$ on $\cC$ (that is the corresponding
parts of their canonical bases given by (\ref{eqBDivision})
are equipotent), symbolically $\varphi|\cB \cong \psi|\cC,$ then
$\varphi$ and $\psi$ are conjugate.

For the sake of simplicity
we prove the converse (in fact a stronger result)
only for involutions we essentially use: for
involutions with $Z(\cB)=\varnothing$ in all
canonical bases $\cB$'s.  We shall call these
involutions {\it soft} involutions.  It is useful that
involutions in $\aut A$ induced by them have a sum of
eigen $\pm$-subgroups equal to $A,$ like
involutions in general linear groups over division rings of
characteristic $\ne 2$.

We shall say that involutions $\varphi,\psi \in \aut F$ {\it have
the same canonical form,} if for all canonical bases
$\cB,\cC$ of $\varphi$ and $\psi,$ respectively, $\varphi|\cB
\cong \psi|\cC.$ Note that {\it a priori} we cannot even claim
that the relation we introduce is reflexive.

\begin{Prop} \label{ConjCriterion}
Let $\varphi \in \aut F$ be a soft involution.  An
involution $\psi \in \aut F$ is conjugate to $\varphi$ if and only if
$\psi$ is soft and $\varphi,\psi$ have the same canonical
form.
\end{Prop}

\begin{pf} Let $A_2 = A/2A$ that is the quotient group
of $A$ by the subgroup of even elements. Natural
homomorphisms $F \to A$ and $A \to A_2,$ gives us a
homomorphism $\mu : \aut F \to \aut{A_2}.$ Clearly, by
\ref{CanForm} the family of all involutions in $\ker
\mu$ coincides with the family of all soft
involutions. Therefore an involution which is
conjugate to a soft involution is soft too.

\begin{Lem} \label{TheyMoveToInverse}
Let $\varphi$ be a soft involution with a canonical basis $\cB.$

{\rm (i)} Suppose $a$ is an element in $F$ such
that $\varphi a=a\inv.$ Then $a=\varphi(w)w\inv$ or $a=\varphi(w) x
w\inv$ for some $x \in X=X(\cB)$ and $w \in F.$

{\rm (ii)} Suppose $C$ is a maximal subgroup of $F$ on which
$\varphi$ acts as conjugation by $x \in X(\cB)${\rm:}
$$
C=\{c \in F : \varphi(c)=xcx\inv\}.
$$
Then $C=C^\varphi_x=\str{Y_x}.$
\end{Lem}

\begin{pf} (i) By \ref{CanForm} we have
\begin{equation} \label{eqDecomp}
F = \Fix(\f) * {\prod_{x \in X}}^\ast\, H_x,
\end{equation}
where each factor is $\f$-invariant. Then
$a=a_1 \ldots a_n,$ where every $a_i,$ $i=1,\ldots,n$
is an element of a free factor in expansion
(\ref{eqDecomp}), and  $a_i$ and $a_{i+1}$ lie in
different factors for every $i=1,\ldots,n-1$ (that
is the sequence $a_1,a_2,\ldots,a_n$ is reduced). Hence if
$\f(a)=a\inv,$ or $\f(a_1)\ldots \f(a_n)a_1 \ldots
a_n=1$ then
$$
\f(a_n)a_1=1,\quad \f(a_{n-1})a_2=1, \quad \ldots
$$
Therefore $a=\f(w_0)w_0\inv$ or $a=\f(w_0)\, v\, w_0\inv,$
where $v \in H_x$ for some $x \in X$ and $\f(v)=v\inv.$

So let us prove, using induction on length of
a word $v$ in the basis $\cB,$ that $v=\f(w_1)w_1\inv$ or
$v=\f(w_1)x w_1\inv.$ The only words $v$ of length one
in $H_x$ with $\f v=v\inv$ are $x$ and $x\inv=\f(x)x x\inv.$

An arbitrary element $v \in H_x$ can be written
in the form
\begin{equation} \label{eqZ_In_H_x}
v=x^{k_1} y_1 x^{k_2} y_2 \ldots x^{k_m} y_m,
\end{equation}
where $y_i \in C_x,$ the elements $x^{k_1}$
and $y_m$ could be equal to 1, but any other
element is non-trivial. Since $\f$ acts
on $C_x$ as conjugation by $x$ we have
\begin{equation}  \label{eqPhiOnZ}
\f(v) =x^{-k_1+1} y_1 x^{-k_2} y_2 \ldots x^{-k_m} y_m x\inv.
\end{equation}
Suppose that $\f(v)v=1.$ We then have
$$
x^{-k_1+1} y_1 x^{-k_2} y_2 \ldots x^{-k_m} y_m x\inv x^{k_1} y_1 x^{k_2} y_2 
\ldots x^{k_m} y_m=1.
$$
Let first $y_m \ne 1.$ Then $k_1=1$ and $y_m=y_1\inv.$
Hence
$$
v=xy_1x\inv (x^{k_2+1} y_2 \ldots x^{k_m}) y_1\inv=\f(y_1) t y_1\inv.
$$
It is easy to see that $\f(t)=t\inv$ and length of $t$ is less
than length of $v.$ In the case when $y_m=1$ we
have $k_m \ne 0$ and $k_1=k_m+1.$ Therefore,
$$
v=x^{k_m} (x y_1 \ldots y_{m-1}) x^{k_m} =\f(x^{-k_m}) t x^{k_m},
$$
and we again have that $\f(t)=t\inv$ and $|t|<|v|.$

(ii) Let $\f c=xcx\inv$ and $c=c_1 c_2\ldots c_n,$ where
$c_i$ are elements in free factors from (\ref{eqDecomp})
and the sequence $c_1,c_2,\ldots,c_n$ is reduced. Suppose
that $n \ge 2.$ Due to the $\f$-invariance of our
free factors, the sequence $\f(c_1),\f(c_2),\ldots,\f(c_n)$
is also reduced and must represent the same element as the
sequence $x,c_1,c_2,\ldots,c_n,x\inv.$ It is easy
to see that it is possible if both $c_1,c_n$ lie
in $H_x.$ In particular, $n \ge 3.$ It implies
that $\f(c_1)=xc_1$ and $\f(c_n)=c_n x\inv.$
It easily follows from (\ref{eqZ_In_H_x}) and
(\ref{eqPhiOnZ}) there is no
$v \in H_x$ with $\f(v)=xv;$ it of course means
that for every $v \in H_x$ the equality
$\f(v)=vx\inv$ is also impossible, since
it is equivalent to $\f(v\inv)=x v\inv.$

Thus, if $\f(c)=xcx\inv,$ then $c \in H_x.$
By applying formulae (\ref{eqZ_In_H_x}) and
(\ref{eqPhiOnZ}), one can readily conclude
that $c$ must be in $C_x.$
\end{pf}

{\sc Remarks.} (a) Note that $a=\f(w)w\inv$ cannot be
a {\it primitive} element (i.e. a member of a basis of $F$) since
$\av a$ is an {\it even} element of $A.$ Indeed, it follows
from (\ref{eqDecomp}) that
$$
\avst w =\av{w(U)} + \av{w(X)} +\av{w(Y)}.
$$
where $w(U) \in \Fix(\f),$ $w(X)$ is an element in the subgroup
generated by $X,$ and $w(Y)$ is an element in the subgroup
generated by the set $Y=\bigcup_{x \in X} Y_x.$
Hence
\begin{align*}
\avst a &= \av{\f(w)w\inv}=\av{\f(w)}-\avst w \\
      &= (\av{w(U)} - \av{w(X)} +\av{w(Y)}) - (\av{w(U)} + \av{w(X)} +\av{w(Y)}) 
\\
      &= -2 \av{w(X)}.
\end{align*}

(b) Using a similar argument, we see that
if $\f(w_1) x_1 w_1\inv=\f(w_2) x_2 w_2 \inv,$ where $x_1,x_2 \in X,$
then $x_1=x_2$ (the element $\av x_1-\av x_2$ is even
if and only if $x_1=x_2$).

Suppose now that soft involution $\psi$ is a conjugate of $\f$:
$\psi=\s\inv \f \s.$ Let
$$
\cB'=U' \cup X'  \cup \bigcup_{x' \in X'} Y'_{x'}
$$
be a canonical basis for $\psi.$ Fixed point subgroups
of $\f$ and $\psi$ are clearly isomorphic. Thus,
$|U'|=|U|.$

If $\psi\, x'=x'{}\inv,$ where $x' \in X',$ then $\f(\s
x')=(\s x')\inv.$ By \ref{TheyMoveToInverse} (i) and the
above remarks there is a unique $x \in X$ such that
\begin{equation}  \label{eqs(x')=...}
\s x'=\f(w) x w\inv
\end{equation}
The mapping $x' \mapsto x$ determined in such a way is
injective, because otherwise we can find two distinct elements
in a basis of a free abelian group whose difference is
even.

Hence, $|X'| \le
|X|$ and by symmetry $|X'|=|X|.$

We claim now that
$$
\s C'_{x'} = wC_x w\inv.
$$
It will complete the proof, because in this
case $|Y'_{x'}|=|Y_x|.$

Let $y' \in C'_{x'}$ and $b=\s y'.$ Since $\psi y' =x' y' x'{}\inv,$
then $\f b= (\s x') b (\s x')\inv.$ Therefore we have
$$
\f(w\inv b w) = \f(w\inv) \f(b) \f(w)=x(w\inv b w) x\inv.
$$
Hence by \ref{TheyMoveToInverse} (ii) $w\inv b w \in C_x.$

The equation (\ref{eqs(x')=...}) can be rewritten as
follows
$$
\s\inv x = \psi(\s\inv w\inv) x' \s\inv w.
$$
Therefore
$$
\s\inv C_x \subseteq \s\inv w\inv C'_{x'} \s\inv w,
$$
or
$$
C_x \subseteq w\inv \s C'_{x'} w,
$$
and the result follows.
\end{pf}

\section{\it Anti-commutative conjugacy classes} \label{acc}

The key roles in the group-theoretic characterization
of conjugations in $\aut F$ will be played by two conjugacy classes
of involutions. First class consists of involutions
with a canonical form
\begin{alignat*}2
\f x &=x \inv, & \quad & \\
\f y &= x y x\inv, & y &\in Y
\end{alignat*}
that is $\cB=\{x\} \cup Y$ is a basis of $F,$ canonical for $\f,$
$U(\cB)=Z(\cB)=\varnothing,$
$X(\cB)$ is a singleton set $\{x\},$ $Y_x=Y$ (any canonical
form reproduced below is interpreted in a similar way).
We shall call these involutions {\it quasi-conjugations.}

An arbitrary element $\f \in \aut F$ in the second
class has the following canonical form:
$$
\f x =x\inv,\quad x \in X,
$$
that is there is a basis of $F$ such that $\f$ inverts
all its elements. We shall use for these involutions
the term {\it symmetries.}

Every symmetry induces in $\aut A$ the automorphism
$-\id_A,$ and hence the product of two symmetries
induces $\id_A.$ Therefore by \ref{CanForm} two
symmetries commute if and only if they are coincide.
Thus, the conjugacy class of all symmetries is, say,
{\it anti-commutative}, since its elements are
pairwise non-commuting. In the next section we shall
prove that the class of all quasi-conjugation is also
anti-commutative.

In order to characterize conjugations we shall use
anti-commutative conjugacy classes of involutions, but we need not an
exact determination of {\it all} such conjugacy classes: it
suffices to know that they lie in some
`easy-to-define' family. In the following proposition
we formulate and prove a necessary condition of being
an anti-commutative conjugacy class, but do not prove its
sufficiency.

\begin{Prop} \label{ACCs}
Let $\f$ be an involution in an anti-commutative
conjugacy class. Then either $\f$ has a canonical
form such that
\begin{alignat}2 \label{eqBeads}
\f u &=u,  &\quad u &\in U,\\
\f x &=x\inv,  & x & \in X,  \nonumber \\
\f y &=x y x\inv,  & y & \in Y_x, \nonumber
\end{alignat}
where $|X| \ge 2,$ all the sets $Y_x,$ $x \in X$ have
the same finite power $n$ and $|U| < n+1,$ or $\f$ has
a canonical form such that
\begin{alignat}2 \label{eqSnakes}
\f u &=u,  &\quad u &\in U,\\
\f x &=x\inv,      &\quad &  \nonumber \\
\f y &=x y x\inv,  & y & \in Y, \nonumber
\end{alignat}
where the cardinal $|U|$ is finite and less than $|Y|+1.$
\end{Prop}

In other words, in the terminology introduced in
Section \ref{invs}, a canonical basis $\cB$ for an
involution in an anti-commutative conjugacy class
either contains exactly one block and the power of the
fixed part of $\cB$ is less than the size of this
block, or all blocks of $\cB$ have the same finite
size and the power of the fixed part is less than the
size of any block. Clearly, symmetries have the form
\eqref{eqBeads} (all blocks of their canonical bases
have the size one), and quasi-conjugations have the
form \eqref{eqSnakes} (the size of the unique block
is equal to $\rank F$).

\begin{pf}
Show first that every anti-commutative class of
involutions consists only of soft involutions. Indeed, let
involution $\f$ be an involution whose canonical basis
$\cB$ has non-empty `permutational' part $Z(\cB).$
Suppose $\f$ takes $z \in Z(\cB)$ to $z'.$
Consider an involution $\psi$ which acts on $\cB
\setminus \{z,z'\}$ exactly as $\f$ does, but taking
$z$ to $z'{}\inv.$ Clearly, $\psi$ is conjugate to
$\f$ and commutes with $\f.$

The following example demonstrates why the size of the fixed part of
a canonical basis must necessarily be less than the size of each block:
$$
\begin{cases}
\f x=x\inv,\\
\f y=x y x \inv,\\
\f x_1=x_1,\\
\f y_1=y_1,\\
\f u=u
\end{cases}
\quad
\begin{cases}
\psi x =x,\\
\psi y =y,\\
\psi x_1 =x_1\inv,\\
\psi y_1 =x_1 y_1 x_1 \inv,\\
\psi u =u.
\end{cases}
$$

Let us make a technical remark. Involutions
\begin{equation}  \label{eqTwoMinuses}
\begin{cases}
\f x=x \inv,  \\
\f y=x y x\inv,\, y \in Y
\end{cases}
\mbox{   and   }
\begin{cases}
\psi x=x \inv, \\
\psi y=x\inv y x\inv,\, y \in Y
\end{cases}
\end{equation}
are conjugate: the second one acts `canonically' on
the set $\{x, x\inv y\ : y \in Y\}$: $\psi(x\inv y)=
x(x\inv y)x\inv.$

Let now $\f$ be a soft involution with a canonical
basis $\cB$ such that for some distinct $x_1, x_2 \in
X(\cB)$ either $|Y_{x_1}| < |Y_{x_2}|$ or
$|Y_{x_1}|=|Y_{x_2}|$ and the cardinal $|Y_{x_1}|$ is
infinite (thus, there must be no neither a pair of
blocks of different size nor a pair of infinite
blocks). The fact that $\f$ commutes with at least two
its conjugates is a consequence of
(\ref{eqTwoMinuses}) and the following

\begin{Clm}
Let $G$ be a free group with a basis
$$
\{x,a\} \cup B \cup \{c\} \cup D \cup E,
$$
where $|B|=|D|.$ Then

{\rm (i)} Involutions
$$
\begin{cases}
\f x=x\inv,\\
\f a=x\inv a x\inv,\\
\f b=x\inv b x\inv,\\
\f c=c\inv,\\
\f d=c\inv d c\inv,\\
\f e=x\inv e x\inv
\end{cases}
\qquad \text{ and }\qquad
\begin{cases}
\psi x=x\inv,\\
\psi a=a\inv,\\
\psi b=a\inv b a\inv, \, b \in B,\\
\psi c=x\inv c x\inv,\\
\psi d=x\inv d x\inv, \, d \in D,\\
\psi e=x\inv e x\inv, \, e \in E
\end{cases}
$$
are conjugate and commute.

{\rm (ii)} Furthermore,
\begin{align*}
&\rank C^{\f}_x = |\{a \} \cup B \cup E| =|B|+|E|+1,\\
&\rank C^{\f}_c=|D|.
\end{align*}
Thus, $\rank C^{\f}_c \le \rank C^{\f}_x,$ where
equality holds only if both ranks are infinite.
\end{Clm}
\begin{pf}
Easy.
\end{pf}

To complete the proof of the Proposition, we should cut off
involutions with a canonical form such that
\begin{alignat*}2
\f u &=u, & u &\in U,\\
\f x &=x\inv, &&\\
\f y &=xyx\inv, &\quad y &\in Y,
\end{alignat*}
where $|U|$ is infinite. To do this let us
consider involutions $\f,\psi$ which act
on a basis $U \cup \{x\} \cup Y$ of $F$
as follows
$$
\begin{cases}
\f u_0=u_0,\\
\f u=u,\\
\f x=x\inv,\\
\f y=xyx\inv,
\end{cases}
\quad
\begin{cases}
\psi u=u_0\inv,\\
\psi u=u,\,        \qquad     &u \in U \setminus \{u_0\}, \\
\psi x=u_0 x u_0\inv,\\
\psi y=u_0 y u_0\inv,  \quad  &y \in Y,
\end{cases}
$$
where $u_0$ is a fixed element in $U.$ It is readily seen that
$\f$ and $\psi$ are conjugate
and commute.
\end{pf}

\section{\it Quasi-conjugations} \label{q-conjs}

In this section we obtain a (first-order)
characterization of quasi-conjugations in $\aut F.$
First we prove that the conjugacy class of all
quasi-conjugations is anti-commutative. Then we distinguish
this class from other anti-commutative conjugacy classes of involutions:
it is trivial in the case when $\rank F=2$
and more technical in the case when $\rank F > 2.$

\begin{Prop} \label{QsConjs:ACC} The class of all
quasi-conjugation is an anti-commutative conjugacy
class.
\end{Prop}

\begin{pf}
Let $\f$ be a quasi-conjugation with a canonical
basis $\cB=\{x\} \cup Y$:
\begin{alignat*}2
\f x &=x \inv, & \quad & \\
\f y &= x y x\inv, & y &\in Y
\end{alignat*}
Every $\s \in \aut F,$ which commutes
with $\f,$ takes $x$ to a primitive element of the form
$\f(w)xw\inv,$ $w \in F$ (Lemma
\ref{TheyMoveToInverse}).  It turns out that the only
primitive elements of the form $\f(w) x w\inv$ are $x$
and $x\inv.$ This fact is a consequence of the
following result, which we shall use once more later.

\begin{Lem} \label{X_And_X_inv}
Let $\alpha$ be an involution with a canonical form
such that
\begin{alignat*}2
\alpha\, x &=x\inv, &\quad&\\
\alpha\, y &=xyx\inv, &y &\in Y=Y_x,\\
\alpha\, u &=u,       &u &\in U,
\end{alignat*}
and $a\in F$ a primitive element, which $\alpha$
sends to its inverse. Then $a=vx^{\pm1} v\inv,$
where $v \in \str U =\Fix(\alpha).$
\end{Lem}

\begin{pf}
As we observed earlier $a=\alpha(w)xw\inv.$
It is easy to see that $a$ lies in the normal
closure of $x.$ One can use induction on
length of a reduced word $w$ in the basis
$\{x\} \cup Y \cup U.$ We have
\begin{align*}
&\alpha(xw_1)xw_1\inv x\inv=x\inv (\alpha(w_1)xw_1\inv) x\inv,\\
&\alpha(uw_1)xw_1\inv u\inv=u(\alpha(w_1)xw_1\inv)u\inv,\\
&\alpha(yw_1)xw_1\inv y\inv= x \cdot y(x\inv \alpha(w_1)xw_1\inv)y\inv.
\end{align*}

\begin{prop} \label{Magnus'sProp}
\mbox{\rm (\cite{Ma}, \cite[II.5.15]{LSch})}
Let $F$ be a free group and the normal closure of $q
\in F$ consists of a primitive element $p.$ Then $q$
is conjugate to $p$ or $p\inv.$
\end{prop}

Hence $a=bx^\eps b\inv,$ where $\eps=\pm 1.$ Since
$\alpha(a)=a\inv,$ we have
$$
\alpha(b) x^{-\eps} \alpha(b\inv)=b x^{-\eps} b\inv.
$$
It follows that $x^{-\eps}$ and $b\inv \alpha(b)$
commute. Therefore both these elements lie
in a cyclic subgroup of $F$ (\cite[I.2.17]{LSch}).
It must be the subgroup $\str x,$ because $x$ is primitive.
Hence $\alpha(b)=b x^k.$ If $k$ is even, say
$k=2m,$ then $\alpha(b x^m)=bx^m$ and $v=bx^m  \in \str U.$
In the case when $k$ is odd, we have  there is
$z \in F$ such that $\alpha(z)=zx.$ One can easily check
that this is impossible (e.g. using `abelian' arguments
as in Section \ref{invs}).
\end{pf}

\begin{Prop} \label{CenOfASnake}
The centralizer of a quasi-conjugation $\f$ with
a canonical basis $\{x\} \cup Y$ consists of automorphisms $\s \in \aut F$
of the form
\begin{alignat}2  \label{eqCenOfASnake1}
\s x &=x,  &\quad & \\
\s y &= \theta(y), & y &\in Y,  \nonumber
\end{alignat}
and of the form
\begin{alignat}2 \label{eqCenOfASnake2}
\s x &=x\inv,  &\quad & \\
\s y &=x \theta(y) x\inv, & y &\in Y, \nonumber
\end{alignat}
where $\theta \in \aut{C_x}.$
\end{Prop}

\begin{pf}
The Proposition easily follows from Lemma \ref{X_And_X_inv}
and one more

\begin{Lem} \label{A*B}
\mbox{\rm \cite[p. 101]{Cohen}}
Let $F=G*H$ be a free factorization of a free group
$F.$ Suppose that $\alpha \in \aut F$ and $\alpha|_{G}$ is an
endomorphism of $G.$ Then $\alpha|_{G} \in \aut G$ if
either $\rank G$ or $\rank H$ is finite.
\end{Lem}

By \ref{X_And_X_inv} if $\s \in \Cen(\f),$ then either
$\s x=x$ or $\s x=x\inv.$ Assume that $\s$ fixes $x.$
Then, for each $y \in C_x$
$$
\f\s y=\s x\, \s y\, \s x\inv= x\, \s y\, x\inv.
$$
Therefore $\s y \in C_x.$

If $\s x=x\inv,$ then $\f\s x=x.$ It follows that $\f\s$ has
the form (\ref{eqCenOfASnake1}), and hence
$\s$ must have the form (\ref{eqCenOfASnake2}).
\end{pf}

We can prove now that $\f$ is the unique
quasi-conjugation in its centralizer, or, in other
words, the conjugacy class of $\f$ is anti-commutative.
Indeed, there are no
quasi-conjugations in the family of automorphisms
of the form (\ref{eqCenOfASnake1}), because every quasi-conjugation
has trivial fixed-point subgroup.
Let $\s$ be a quasi-conjugation of the form
(\ref{eqCenOfASnake2}). Hence $\theta^2=\id,$ and
$\theta$  is either the identity automorphism of $C_x$ or a
soft involution. Assume that $\theta$ is an
involution. By Theorem \ref{CanForm} there is a
canonical basis $\cC$ of a free group $C_x$ for
$\theta$ such that
\begin{alignat*}3
&\theta & u &=u,              & u &\in U(\cC),\\
&\theta & c &=c\inv,          & c &\in X(\cC),\\
&\theta & d &=c d c\inv,\quad & d &\in Y_c(\cC).
\end{alignat*}
Therefore
\begin{alignat*}3
&\s  x    &&=x\inv,                                  \\
&\s  u    &&=x u x\inv,              & u &\in U(\cC),\\
&\s  (xc) &&=(xc)\inv,               & c &\in X(\cC),\\
&\s  d    &&=(xc) d (xc)\inv,\quad   & d &\in Y_c(\cC).
\end{alignat*}
We obtain a canonical basis for $\s.$ This basis
contains at least two blocks, and hence by Proposition
\ref{ConjCriterion} $\s$ cannot be a quasi-conjugation.
So $\theta$ must be the identity automorphism, or, equivalently,
$\s=\f$ as desired. This completes the proof. \end{pf}

Now there are no further problems in a characterization
of quasi-conjugations in automorphism groups of two-generator
free groups.

\begin{Prop} \label{rank(F)=2}
Suppose that $\rank F=2.$ Then
the class of all quasi-con\-ju\-ga\-tions is the
unique anti-commutative conjugacy class $K$ of
involutions such that elements in $K$ are not squares.
\end{Prop}

\begin{pf}
By Propositions \ref{ACCs} and \ref{QsConjs:ACC} the
only anti-commutative conjugacy classes of involutions in
$\aut F$ are the class of all quasi-conjugations and the
class of all symmetries.  Every symmetry is a square:
$$
\begin{cases}
\psi x_1=x_2\inv, \\
\psi x_2=x_1
\end{cases}
\Rightarrow
\begin{cases}
\psi^2 x_1=x_1\inv,\\
\psi^2 x_2=x_2\inv.
\end{cases}
$$
On the other hand, a quasi-conjugation $\f$ with a canonical
form
\begin{align*}
\f x &=x\inv,\\
\f y &=xyx\inv
\end{align*}
induces in $\aut A$ an automorphism with the determinant which is equal
to $-1.$ Hence $\av\f$ cannot be a square in $\aut A.$ It implies
that $\f$ is not a square in $\aut F.$
\end{pf}

Let us consider now more serious case.

\begin{Th} \label{DefOfQsConjs}
Let $\rank F > 2.$ Then the class of all quasi-conjugations is the
unique anti-commutative conjugacy class $K$ of
involutions such that for every anti-commutative
conjugacy class $K'$ of involutions, all involutions in $KK'$ are
conjugate.
\end{Th}

\begin{pf} Let $\f$ be a quasi-conjugation, $K'$
an anti-commutative conjugacy class of involutions, and
$\f \not\in K'.$ It suffices to prove that if
$\psi,\psi' \in K'$ both commute with $\f,$ then
$\f\psi$ and $\f\psi'$ are conjugate.

Suppose first that fixed point
subgroups of $\psi$ and $\psi'$ are trivial. Then both
$\psi$ and $\psi'$ have the form
\eqref{eqCenOfASnake2}:
\begin{alignat*}2
\psi x &=x\inv,                  & \psi' x &=x\inv \\
\psi y &=x \theta(y) x\inv,\quad & \psi' y &=x\theta'(y) x\inv, \quad y \in Y,
\end{alignat*}
where $\theta$ and $\theta'$ are in $\aut{C_x}.$ It is easy to see that
$\theta$ and $\theta'$ are soft involutions. Let $\cC$ be a canonical
basis for $\theta.$ As we observed above, the set
\begin{equation}
(\{x\} \cup U(\cC)) \cup \bigcup_{c \in X(\cC)} (\{xc\} \cup Y_c(\cC))
\end{equation}
is a canonical basis for $\psi,$ and (\theequation) is a partition
of this basis into blocks (there are at least two blocks,
because $\psi$ cannot be conjugate to $\f$). Since $\psi$ lies
in an anti-commutative conjugacy class,
all these blocks have the same size,
and hence for all $c \in X(\cC)$
$$
|U(\cC)|=|Y_c(\cC)|.
$$
By applying a similar argument to $\psi',$ we see
that $\theta$ and $\theta'$ are conjugate in $\aut {C_x}.$
Hence $\f \psi$ and $\f \psi'$ are conjugate in $\aut F$:
both these automorphisms fix $x$ and their restrictions
on $C_x$ ($\theta$ and $\theta'$) are conjugate.

Assume that $\psi$ and $\psi'$ have non-trivial fixed
point subgroups. Both $\psi$ and $\psi'$ preserve the
subgroup $C_x$ and fix the element $x.$ It means
(Propositions \ref{CanForm}, \ref{ConjCriterion}) that restrictions of
$\psi$ and $\psi'$ on $C_x,$ say, $\theta$ and
$\theta',$ respectively are conjugate in $\aut{C_x}$:
$\theta'=\pi\inv \theta \pi,$ where $\pi \in \aut{C_x}.$
The automorphism $\pi$ can be extended to an element
$\s \in \aut F$ such that $\s x=x.$ By \ref{CenOfASnake}
$\s$ commutes with $\f,$ and hence
$$
\s\inv(\f \psi)\s=\f \psi'.
$$

Let us prove the converse. We start with the following

\begin{Clm} \label{PlaysWithSymmetries}
Assume that an involution $\f \in \aut F$
is an element of an anti-commutative conjugacy class, and $\f$ is neither
quasi-conjugation, nor symmetry. Then there exist
symmetries $\psi$ and $\psi'$ such that $\f\psi$ and
$\f\psi'$ are non-conjugate involutions.
\end{Clm}

\begin{pf} (a) A natural way to define an action of a symmetry commuting
with $\f$ on a block of a canonical basis for $\f$
is given by Proposition \ref{CenOfASnake}:
\begin{alignat*}2
\f x &=x\inv,           & \psi x &=x\inv,\\
\f y &=x y x\inv \quad  & \psi y &=x y\inv x\inv,\quad (\psi(xy)=(xy)\inv) \quad 
y \in Y_x.
\end{alignat*}
Clearly, the product of $\f$ and $\psi$ fixes $x$ and inverts
each element in $Y_x.$

(b) A way to define an action of a symmetry on the fixed
part of a canonical basis is obvious: a symmetry should invert
each element.

If a symmetry $\psi$ acts on a canonical
basis for $\f$ as it is defined in (a) and (b), then
any canonical basis of the product $\f\psi$ has non-empty
fixed part and each block of this basis has the size one.

(c) Assume now that a canonical basis $\cB$ for $\f$ contains
at least two blocks, say,
\begin{equation}
(\{x\} \cup (\{y\} \cup B)) \text{ and } (\{z\} \cup (\{t\} \cup C))
\end{equation}
Hence
\begin{alignat*}2
\f x &=x\inv,\\
\f y &=x y x\inv,\\
\f a &=x b x\inv, \quad & b &\in B,\\
\f z &=z\inv,\\
\f t &=ztz\inv,\\
\f c &=z c z\inv, \quad & c &\in C
\end{alignat*}
Then we can choose a symmetry $\psi$
such that any canonical basis for $\f\psi$ has a block of the
size two. Let $\cC$ denote the union of blocks in (\theequation).
Suppose that $\psi$ acts on $\cB \setminus \cC$ as it is defined
in (a) and (b), and define its action on the subgroup generated
by $\cC$ as follows
\begin{alignat*}4
&\psi && (tx) &&=(tx)\inv,\\
&\psi && y    &&=(tx)y\inv (tx)\inv,\\
&\psi && b    &&= (tx) b \inv (tx)\inv, \quad & b &\in B,\\
&\psi && z    &&=z\inv,\\
&\psi && t    &&=zt\inv z\inv,\\
&\psi && c    &&=z c\inv z\inv, \quad & c &\in C
\end{alignat*}
Let $\s$ denote the automorphism $\f\psi.$  We have
$$
\s(tx)=\f(x\inv t\inv)=x zt\inv z\inv.
$$
Since $\s z=z$ and $\s t=t\inv,$ then
\begin{equation} \label{eq2-Block}
t\inv \s(x) z =(xz) t\inv \iff \s(xz)=t (xz) t\inv.
\end{equation}
Furthermore,
\begin{equation}
\s(y)= \f(t x y\inv x\inv t\inv)=\f(t)\f(x y\inv x\inv) \f(t) =
ztz\inv y\inv zt\inv z\inv.
\end{equation}
Let $y'$ denote the element $ztz\inv y.$ It easily follows
from (\theequation) then $\s$ inverts $y'.$
The same is true for all elements $b'=ztz\inv b,$ where $b\in B.$
Summing up, we conclude that $\s$ is an involution and
the basis
$$
(\cB \setminus \cC) \cup (\{t\} \cup \{xz\}) \cup \{y'\} \cup
\{b' : b \in B\} \cup \{z\} \cup C
$$
is a canonical basis for $\s.$ The formulae in (\ref{eq2-Block}) demonstrate
that this basis contains a block of the size two.

(d) Any canonical basis for $\f$ has exactly one block and
non-empty fixed part:
\begin{alignat*}2
\f u &=u, &\quad u &\in U,\\
\f x &=x\inv, && \\
\f y &=xyx\inv, & y &\in Y
\end{alignat*}
One can easily find a symmetry $\psi$ such that the product
$\f\psi$ will be non-conjugate to each product of $\f$ with a symmetry
obtained in a natural way, using (a) and (b):
\begin{alignat*}4
&\psi && u_0  &&=u_0\inv, &&\\
&\psi && u    &&=u\inv, &  u &\in U \setminus \{u_0\},\\
&\psi && x    &&=u_0 x\inv u_0\inv, &&\\
&\psi && y    &&=u_0 xy\inv x\inv u_0\inv, &\quad y &\in Y,
\end{alignat*}
where $u_0 \in U.$ The reason is the same as in the previous
point (c): any canonical basis for $\f\psi$ contains
a block of the size two ($\{u_0\} \cup \{x\}$ in this
example).
\end{pf}

To complete the proof of the Theorem, we have to find
for an arbitrary symmetry $\f$ two involutions
$\psi,\psi'$ from an anti-commutative conjugacy class such that
their products with $\f$ are non-conjugate
involutions. Clearly, the problem is to add one more
conjugacy class to the family of definitely known at
this moment anti-commutative conjugacy classes of involutions
(quasi-conjugations and symmetries). Having such a
class, one can rework in an obvious way examples in
the proof of the latter Claim, and hence prove the
desired result on symmetries. It will complete the proof
of Theorem \ref{DefOfQsConjs}.

\begin{Clm} \label{OneMoreACC}
Let $\rank F > 2.$ The conjugacy class of an involution
$\psi$ with a canonical form such that
\begin{alignat*}2
\psi x &=x\inv, &\quad&\\
\psi y &=xyx\inv, &y &\in Y=Y_x,\\
\psi u &=u,
\end{alignat*}
where $Y \ne \varnothing$ is anti-commutative {\rm(}the fixed point
subgroup of $\psi$ is of rank one{\rm)}.
\end{Clm}

\begin{pf}
Suppose that $\s \in \Cen(\psi).$ The only primitive
elements which are fixed by $\psi$ are $u$ and
$u\inv.$ Therefore $\s u=u^{\pm 1}.$ By Lemma
\ref{X_And_X_inv} $\s$ must take $x$ to an element of
the form $u^k x^{\pm 1} u^{-k},$ where $k \in \Z.$ To
calculate the image $\s y$ of an element $y \in Y$ one
may apply the arguments used in the proof of
Proposition \ref{CenOfASnake}. We then have that
$\s$ has the form
\begin{alignat}2  \label{eqCenOfA1-Snake1}
\s u &=u^{\varepsilon},  &\quad & \\
\s x &=u^k x u^{-k},  &\quad & \nonumber \\
\s y &= u^k\theta(y) u^{-k}, & y &\in Y,  \nonumber
\end{alignat}
or the form
\begin{alignat}2 \label{eqCenOfA1-Snake2}
\s u &=u^{\eta},  &\quad & \\
\s x &=u^m x\inv u^{-m},  &\quad & \nonumber\\
\s y &=u^m x \theta(y) x\inv u^{-m}, & y &\in Y, \nonumber
\end{alignat}
where $\varepsilon,\eta=\pm1,$ $k,m$ are integers, and
$\theta \in \aut{C_x}.$ Note that if ($\varepsilon=1$ and $k\ne 0$)
or ($\eta=1$ and $m \ne 0$), then $\s$ has infinite order.

Let now $\s \in \Cen(\psi)$ be a conjugate of
$\psi.$ We should prove that $\s=\psi.$

Suppose first that $\s$ has the form \eqref{eqCenOfA1-Snake1}.
When $\varepsilon=1$ and $k=0,$ the subgroup $\Fix(\s)$
is of rank at least two and $\s$ cannot be conjugate to $\psi,$ because
$\rank \Fix(\psi)=1.$ Let $\varepsilon=-1;$ then $\theta^2=\id.$
As we observed above it follows from the assumption $\theta \ne \id$ that
a canonical basis for $\s$ contains at least two blocks.
It is impossible. Thus, $\theta=\id.$ When $k$ is even, say $k=2l,$
the subgroup $\Fix(\s)$ has rank at least two:
\begin{alignat*}2
\s u &=u\inv,  &\quad & \\
\s (u^l x u^{-l}) &=u^l x u^{-l},  &\quad &  \\
\s (u^l y u^{-l}) &= u^l y u^{-l}, & y &\in Y,
\end{alignat*}
Similarly, we see that in the case when $k$ is
odd, $\s$ is a quasi-conjugation.

Let us now try to find a conjugate $\s$ of $\psi$ in the family
of automorphisms of the form \eqref{eqCenOfA1-Snake2}. An
involution of this form with $\eta=-1$ is not conjugate to
$\psi$: any canonical basis for such an involution contains
more than one block. Thus, $\eta=1,$ and we have that $m=0$
and $\theta=\id.$ Therefore $\s=\psi.$
\end{pf}

The proof of Theorem \ref{DefOfQsConjs} is now complete.
\end{pf}

Proposition \ref{rank(F)=2} and Theorem \ref{DefOfQsConjs}
can be summarized in model-theoretic terms as follows.

\begin{Th}
The set of all quasi-conjugations is first-order parameter-free
definable in $\aut F.$
\end{Th}

\begin{pf}
All the hypotheses in Proposition \ref{rank(F)=2} and
Theorem \ref{DefOfQsConjs} are in fact first-order.
For instance, the following formula says that the conjugacy
class of an involution $v$ is anti-commutative:
$$
\operatorname{ACC}(v)=(v \ne 1 \land v^2=1) \land (\forall u)((v v^u)^2=1 
\rightarrow v=v^u),
$$
where $v^u=u v u\inv.$ There are no difficulties in conversion of
other hypotheses into first-order formulae.
\end{pf}

{\sc Remarks.} (a) In the case when $F$ has infinite
rank, a group-theoretic characterization of
quasi-conjugations can be obtained in an easier
way. Indeed, involutions in anti-commutative conjugacy
classes of the form
\eqref{eqBeads} (infinitely many finite blocks of the
same size in any canonical basis) are in this case
squares in $\aut F$:
$$
\begin{cases}
\s x =a\inv,\\
\s y =a b a\inv,\\
\s a =x,\\
\s b =y,
\end{cases}
\Rightarrow
\begin{cases}
\s^2 x =x\inv,\\
\s^2 y =x y x\inv, \, &y \in Y    \\
\s^2 a =a\inv,                       \\
\s^2 b   =a b a\inv       &b \in B \quad (|B|=|Y|).
\end{cases}
$$
On the other hand, it is easy to see that involutions
of the form \eqref{eqSnakes} (exactly one block in any
canonical basis) are not squares. In particular, the
condition of being a square distinguishes symmetries
from quasi-conjugations. Therefore in order to
characterize quasi-conjugation in $\aut F$ one can use
the following

\begin{th}
Let $F$ has infinite rank. Then the class of all
quasi-conjugations is the unique anti-commutative
conjugacy class $K$ of involutions such that its
elements are not squares and for every anti-commutative
conjugacy class $K'$ of involutions, whose elements are
squares, all involutions in $KK'$ are conjugate.
\end{th}

A proof of the latter Theorem may follow the plan of
the proof of Theorem \ref{DefOfQsConjs}, but there is
no need to consider the point (c) in the proof of Claim
\ref{PlaysWithSymmetries} and Claim \ref{OneMoreACC}.

(b) One can easily obtain a {\it uniform} first-order
characterization of quasi-con\-ju\-ga\-tions. Indeed, the
result of S.~Meskin \cite{Me},\cite[I.4.6]{LSch} states that if $\rank F=2,$ 
then
the group $\aut F$ has exactly
four conjugacy classes of involutions.  The converse
is a consequence of Proposition \ref{ConjCriterion},
since if $\rank F>2,$ then $\aut F$ has at least
six conjugacy classes of soft involutions (the number
of conjugacy classes of soft involutions in the
group $\aut{F_3},$ where $F_3$ is a three-generator
free group). Suppose that first-order formulae $\operatorname{QC}_0(v)$
and $\operatorname{QC}_1(v)$ define quasi-conjugations
in $\aut F$ in the case when $\rank F=2$ and $\rank F>2,$
respectively. Then a first-order formula
$$
(\operatorname{QC}_0(v) \land \chi) \lor (\operatorname{QC}_1(v) \land \lnot 
\chi),
$$
where a closed first-order formula $\chi$ says about
four conjugacy classes of involutions, defines quasi-conjugations in
the automorphism group of an arbitrary non-abelian
free group.

\section{\it Conjugations} \label{conjs}

The following theorem is our key result in the proof
of completeness of $\aut F.$

\begin{Th} \label{DefOfConjsByPPE}
The set of all conjugations by powers of primitive elements
is first-order parameter-free definable in $\aut F.$
\end{Th}

\begin{pf}
As in the proof of definability of quasi-conjugations
we consider two cases: $\rank F=2$ and $\rank F > 2.$

{\it I. $F$ is of rank two.}

Let $\f$ be a quasi-conjugation with a canonical
basis $\cB=\{x\} \cup \{y\}.$ By Proposition \ref{CenOfASnake}
the centralizer of $\f$ consists of four elements.
Non-trivial ones are involutions: $\f,$ a symmetry,
and an involution $\psi$ which fixes $x$ and inverts
$y.$ Clearly, $\psi$ is the unique involution in $\Cen(\f)$
commuting with at least two its conjugates.

\begin{Prop}  \label{Conjs:rank=2}
All conjugations by powers of $x$ are in $\Cen(\psi).$
An element $\s$ in $\Cen(\psi)$ is conjugation by a power of $x$
if and only if $\s$ is not an involution, and can be represented as the product
of two conjugate involutions.
\end{Prop}

\begin{pf}
Let $\s \in \Cen(\psi).$ Since $x$ and $x\inv$ are the only
primitive elements in $\Fix(\psi),$ we have $\s x=x^{\pm 1}.$
We can use then Lemma \ref{X_And_X_inv} (assuming that the size of a block
is equal to one). Therefore $\s y=x^k y^{\pm 1} x^{-k},$
for some $k \in \Z.$ Thus,
\begin{align*}
\s x &=x^\varepsilon,\\
\s y &=x^k y^\eta x^{-k},
\end{align*}
where $\varepsilon,\eta=\pm 1$ (and conversely, every automorphism
of $F$ of the latter form commutes with $\psi$). In the case
when $\varepsilon=-1,$ $\s$ is an involution. The
automorphism $\s$ of $F$ such that
\begin{align*}
\s x &=x,\\
\s y &=x^k y\inv x^{-k}
\end{align*}
induces in $\aut A$ an automorphism  with determinant $-1.$
On the other hand, the product of two conjugate involutions
from $\aut F$ induces in $\aut A$ an automorphism whose determinant
is equal to $1.$ To complete the proof, we should express as the product
of two conjugate involutions an arbitrary conjugation by a power of $x.$
It is easy:
\begin{alignat*}2
\alpha x &=x\inv,                  & \alpha' x &=x\inv \\
\alpha y &=y\inv, \qquad           & \alpha' y &=x^{-k} y\inv x^k.
\end{alignat*}
The product of symmetries $\alpha$ and $\alpha'$ is evidently
conjugation by $x^k.$
\end{pf}

{\it II. $F$ is of rank at least three.}

We also start with a quasi-conjugation $\f.$
Suppose that $\cB=\{x\} \cup Y$ is a canonical basis for
$\f.$ Let $\Pi$ denote the set of all automorphisms of $F$
of the form $\pi=\s\s',$ where $\s$ and $\s'$ are in $\Cen(\f)$
and conjugate. By \ref{CenOfASnake} conjugate automorphisms
$\s,\s'$ in the centralizer of $\f$ either both have
the form \eqref{eqCenOfASnake1} (when their fixed point
subgroups are non-trivial) or have the form \eqref{eqCenOfASnake2}.
Therefore every $\pi \in \Pi$ has the form
\begin{align*}
\pi x &=x\\
\pi y &=\theta(y),\quad y \in Y,
\end{align*}
where $\theta \in \aut{C_x},$ that is $\pi$ fixes $x$ and
preserves the subgroup $C_x.$

\begin{Prop} \label{Conjs:rank=>3}
All conjugations by powers of $x$ are in the centralizer
of the family $\Pi.$ Every member of $\Cen(\Pi)$
is either an involution or conjugation by a power of $x.$
\end{Prop}

\begin{pf} Let
$$
\cC=\{a,b\} \cup C
$$
be a basis of $C_x,$ and $\tau \in \aut F$ an element
of $\Cen(\Pi).$

First we construct $\pi \in \Pi$ such that the
fixed point subgroup of $\pi$ is the subgroup
$\str{x,a}.$ Since $\tau$ must commute with $\pi$ we shall have that
$$
\tau a=w_a(x,a),
$$
where $w_a$ is a reduced word in letters $x$ and
$a.$

To construct $\pi,$ we use the same idea as in the proof of
the previous result:
\begin{alignat*}2
\s x &=x,                   & \s' x &=x, \\
\s a &=a\inv,               & \s' a &=a\inv,\\
\s b &=b\inv,               & \s' b &=a\inv b\inv a,\\
\s c &=c\inv, \qquad        & \s' c &=a\inv c\inv a,\quad c \in C.
\end{alignat*}
The restriction of $\pi=\s\s'$ on $C_x$ is conjugation
by $a.$ Then it is easy to show that the fixed point
subgroup of $\pi$ is $\str{x,a}.$ By Lemma \ref{A*B}
$\tau \str{x,a}=\str{x,a}.$

A similar argument can be applied to an arbitrary
primitive element in $C_x.$ Hence for every primitive
$d \in C_x$
$$
\tau d=w_d(x,d),
$$
and $\tau$ preserves the subgroup $\str{x,d}.$ We then have
$$
\tau \str x=\tau(\str{x,a} \cap \str{x,b})=\str{x,a} \cap \str{x,b}=\str x.
$$
Therefore $\tau x=x^{\pm 1}.$ In particular, the word $w_a(x,a)$
must have explicit occurrences of $a.$

We claim now that the words $w_d,$ where $d=a,b,ab$ have the same structure,
that is any word $w_d(x,d)$ can be obtained from the word
$w_a(x,a)$ by replacing occurrences of $a$ by $d$:
$$
[w_a(x,a)]^a_d=w_d(x,d).
$$
To prove this, it suffices to find in $\Pi$ automorphism of $F$
which takes $a$ to $b$ ($a$ to $ab$).

Let $\s_1$ and $\s'_1$ be involutions in $\Cen(\f)$
such that $\s_1$ and $\s'_1$ both fix
the set $\{x\} \cup C$ pointwise and
\begin{alignat*}2
\s_1 a &=b\inv,          & \s'_1 a &=ab, \\
\s_1 b &=a\inv, \qquad   & \s'_1 b &=b\inv,\\
\end{alignat*}
Clearly, $\s_1$ and $\s'_1$ are conjugate and
$\pi_1=\s_1\s'_1$ sends $a$ to $b.$ Since $\tau$
and $\pi_1$ commute, we have
\begin{align*}
\tau a=w_a(x,a) &\Rightarrow \tau (\pi_1 a) =w_a(\pi_1 x,\pi_1 a) \\
		&\Rightarrow w_b(x,b)=w_a(x,b).
\end{align*}
Thus, there is a reduced word $w$ in letters $x$ and, say, $t$ such that
$$
[w(x,t)]^t_d=w_d(x,d),
$$
where $d=a,b,ab.$ We then have
$$
\tau (ab)=w(x,ab)=\tau(a) \tau (b)=w(x,a)w(x,b),
$$
and hence
\begin{equation}
w(x,ab)=w(x,a) w(x,b).
\end{equation}

Now we show that the word $w(x,t)$ has the form
$x^k t x^{-k},$ where $k \in \Z.$ Assume that $w(x,t)$ has
the (possibly non-reduced) form
such that
$$
x^{k_1} t^{l_1} x^{k_2} t^{l_2} \ldots x^{k_m} t^{l_m},
$$
where $k_1$ or $l_m$ could be equal to zero, whereas any other
exponent is non-trivial. Then by (\theequation)
\begin{equation}
x^{k_1} (ab)^{l_1} x^{k_2} (ab)^{l_2} \ldots x^{k_m} (ab)^{l_m} =
x^{k_1} a^{l_1} x^{k_2} a^{l_2} \ldots x^{k_m} a^{l_m}
            x^{k_1} b^{l_1} x^{k_2} b^{l_2} \ldots x^{k_m} b^{l_m}.
\end{equation}
The latter equality is evidently impossible when $m \ge 2$ and $l_m \ne 0.$
Hence $l_m=0$ and $k_m=-k_1.$ Even after this reduction (\theequation)
fails, if $m \ge 3.$ Therefore
$$
x^{k_1} (ab)^{l_1} x^{-k_1}=x^{k_1} a^{l_1} b^{l_1} x^{-k_1},
$$
and we have
$$
(ab)^{l_1}=a^{l_1} b^{l_1}.
$$
Since $a$ and $b$ are independent, $l_1=1.$

Summing up, we see that $\tau$ acts on $\cB$ as follows
\begin{align}
\tau x &=x^\varepsilon,  \\
\tau a &=x^k a x^{-k}, \nonumber  \\
\tau b &=x^k b x^{-k},   \nonumber\\
\tau c &=x^k c x^{-k}, \quad c \in C. \nonumber
\end{align}
where $\varepsilon=\pm 1.$ In the case when $\varepsilon=-1,$ $\tau$ is an 
involution,
otherwise $\tau$ is conjugation by $x^k.$ Conversely, every
automorphism of $F$ of the form (\theequation) is in $\Cen(\Pi).$
\end{pf}

Using Theorem \ref{DefOfQsConjs}, one can readily
convert the hypotheses in Proposition
\ref{Conjs:rank=2} and Proposition \ref{Conjs:rank=>3}
into first-order formulae.  The proof of Theorem
\ref{DefOfConjsByPPE} is complete.
\end{pf}

{\sc Remark.} In the case when $F$ has {\it infinite}
rank the subgroup of all conjugations is first-order
definable subgroup of $\aut F.$ Indeed, it is easy to
see that every element in infinitely generated free
group can be expressed as the product of two primitive
elements.  Therefore every conjugation in $\aut F$ is
the product of two conjugations by powers of
primitive elements.  This argument of course does
not work for finitely generated free groups.

So all is now in readiness for a proof of the main result of the
paper.

\begin{Th} \label{MainRes}
Let $F$ be a non-abelian free group. Then the group $\aut F$
is complete.
\end{Th}

\begin{pf}
By Theorem \ref{DefOfConjsByPPE} and Claim
\ref{Def=>Char} the subgroup of $\aut F$ generated by
all conjugations by powers of primitive elements, namely the
subgroup $\Inn F$ of all conjugations is a characteristic
subgroup of $\aut F.$ Therefore the group $\aut F$ is complete
(Proposition \ref{Specht?})
\end{pf}

\begin{Th}
The automorphism groups of free groups $F$ and $F'$
are isomorphic if and only if $F \cong F'.$
\end{Th}

\begin{pf}
We can assume that both groups $F,F'$ have ranks at least two.
Any isomorphism from $\aut F$ to $\aut{F'}$ preserves
conjugations, and hence induces an isomorphism
between $F$ and $F'.$
\end{pf}



\begin{thebibliography}{99}


\bibitem{Cohen} R.~Cohen, `Classes of automorphisms of
free group of infinite rank', {\it Trans. Amer. Math. Soc.}, 177
(1973) 99--119.

\bibitem{DFo} J.~Dyer, E.~Formanek, `The automorphism
group of a free group is complete', {\it J.  London Math.
Soc.},  11 (1975) 181--190.


\bibitem{DSc} J.~Dyer, G.~P.~Scott, `Periodic automorphisms
of free groups', {\it Comm. Algebra}0, 3  (1975) 195--201.

\bibitem{Fo}  E.~Formanek, `Characterizing a free group in
its automorphism group', {\it J. Algebra}, 133 (1990)
424--432.

\bibitem{Ho} W. Hodges, `Model Theory', (University
Press, Cambridge, 1993).


\bibitem{Khr}  D.~G.~Khramtsov, `Completeness of groups
of outer automorphisms of free groups', {\it in}
Group-theoretic investigations (Russian),
Akad. Nauk SSSR Ural. Otdel., Sverdlovsk, (1990) 128--143.


\bibitem{LSch} R.~Lyndon, P.~Schupp, `Combinatorial group
theory', (Springer-Verlag, Berlin, etc., 1977).


\bibitem{Ma}  W.~Magnus, `Untersuchungen \"uber einige unendliche
discontinuierliche Gruppen', {\it Math. Ann.}, 105 (1931), 52--74.

\bibitem{Me} S.~Meskin, `Periodic automorphisms of the
two-generator free group', Proc. Conf. Canberra 1973,
(Lecture Notes in Math., 372, 494--498, Berlin, etc.,
Springer).


\bibitem{Spe} W.~Specht, Gruppentheorie (Berlin-G\"ottingen-Heidelberg,
1956).


\bibitem{To} V.~Tolstykh, `Puissance et pl\'enitude:
interpr\'etation des groupes d'automorphismes des
groupes libres', Quatri\`eme Colloque Franco-Touranien
de Th\'eorie des Mod\`eles, Marseille-Lumini 1997,
R\'esum\'es des Conf\'erences, Institut Girard
Desargues, Universit\'e Claude Bernard, Lyon 1, 19.


\end{thebibliography}
\end{document}